\documentclass[11pt]{amsart}
\usepackage[latin1]{inputenc}
\usepackage{amsmath}
\usepackage{amsfonts}
\usepackage{amssymb}

\newtheorem{thm}{Theorem}[section]
\newtheorem{prop}[thm]{Proposition}

\newtheorem{cor}[thm]{Corollary}
\newtheorem{defn}[thm]{Definition}

\begin{document}

\title{Complex Equiangular Cyclic Frames and Erasures}

\author{Deepti Kalra}

\address{Deepti Kalra: Department of Mathematics, University of Houston, 4800
Calhoun Road, Houston, TX 77204-3008 U.S.A.}

\email{deepti@math.uh.edu}

\maketitle

\begin{abstract}

\noindent We derive various interesting properties of complex equiangular cyclic frames for many pairs $(n,k)$ using Gauss sums and number theory. We further use these results to study the random and burst errors of some special cases of complex equiangular cyclic $(n,k)$ frames.\\



\end{abstract}

\section{Introduction}

\noindent Two-uniform frames are of key importance in coding and decoding of vectors. The two-uniform frames have been discussed in \cite{2} where it is proved that when such frames exist, they are optimal for two or more erasures. It is also shown in \cite{2} that a frame is two-uniform if and only if it is equiangular in the terminology of \cite{6}.\\

\noindent It is known that equiangular $(n,k)$-frames, i.e. equiangular frames of $n$ vectors for a $k$-dimensional Hilbert space, can only exist for certain pairs of integers $(n,k)$.  For real Hilbert spaces, necessary and sufficient conditions for the existence of real equiangualr cyclic frames are expressed in terms of the existence of certain types of graphs \cite{2}. It is shown in \cite{9} that the existence of complex equiangular cyclic $(n,k)$-frames depends on the existence of certain difference sets. \\

\noindent Although in \cite{9} the authors show that the necessary and sufficient condition for existence of an equiangular cyclic $(n,k)$-frame is the existence of a corresponding $(n,k,\lambda)$ difference set, the construction of some types of equiangular cyclic $(n,k)$-frames shown in this paper reveal more important properties which are not reflected in \cite{9}. Using the technique in this paper, we can prove results about the error for some equiangular cyclic$(n,k)$-frames. The construction in this paper, along with results for Gauss sums, give the precise information about the error and the corresponding correlation matrix. \\

\noindent In this paper we introduce Gauss sums in the field of frame theory. We observe that the construction of these equiangular cyclic $(n,k)$-frames is a particular kind of Gauss sums and hence we use the theory developed for the same. We also look into the geometry of these frames and show that these frames form a spherical 1-design.\\ 

\section{Basic Concepts and Definitions}

\noindent We begin by the definition of a {\em frame} for a Hilbert space.

\begin{defn}

Let $\mathbb{H}$ be a real or complex Hilbert space and let $F=\lbrace f_i \rbrace_{i \in I}$ be a subset of $\mathbb{H}$, where $I$ is an index set. Then $F$ is called a {\em frame} for $\mathbb{H}$ provided that there are two positive integers $A$,$B$ such that the inequalities 
$$A\Arrowvert x \Arrowvert^2 \leq \sum _{j \in I} \arrowvert \langle x,f_j \rangle \Arrowvert x \Arrowvert^2 \leq B\Arrowvert x \Arrowvert^2$$

\noindent hold for every $x \in \mathbb{H}$. Here $\langle .,. \rangle $ denotes the inner product of two vectors in $\mathbb{H}$.

\end{defn}

\noindent If $A=B=1$, then $\lbrace f_i \rbrace_{i \in I}$ is called a {\em Parseval} frame or {\em Unit Normalized Tight} frame or UNTF. A frame is called {\em uniform} or {\em equal-norm} provided there is a constant $c$ such that $\Arrowvert f_i \Arrowvert = c$ for each $i \in I$.\\

\noindent Let $\mathbb{F}$ be a field of real or complex numbers. Let $\mathbb{F}(n,k)$ be the collection of all Parseval frames for a $k$-dimensional Hilbert space $\mathbb{F}^{k}$ consisting of $n$ vectors. Such frames are called $(n,k)$-frames. The ratio of $n/k$ is called the $redundancy$ ratio  of the $(n,k)$-frame.\\

\noindent It is known that a Parseval frame satisfies the Parseval identity, 
$$ x=\sum _{l\in I} \langle x,f_{l} \rangle f_{l} \ \ \ \forall \ x \in \mathbb{H}.$$

\noindent Let us now look at the case of losing $m$ coefficients, {\em i.e.} the case of {\em m-erasures}. We define the error operator $ E_{i_1,\dots,i_m} $ as 
$$ E_{i_1,\dots,i_m} (x)= x-\sum_{l \neq i_1, \dots,i_m} \langle x,f_{l} \rangle f_{l}=\sum_{j=1}^{m}\langle x,f_{i_j} \rangle f_{i_j}.$$
\\

\noindent The norm of this error operator is given by operator norm of the $m \times m $ correlation matrix

$$
\left[
\begin{array}{cccc}
\langle f_{i_1}, f_{i_1} \rangle & \langle f_{i_2}, f_{i_1} \rangle & \dots & \langle f_{i_m}, f_{i_1} \rangle \\
\langle f_{i_1}, f_{i_2} \rangle  & \langle f_{i_2}, f_{i_2} \rangle &  \ldots & \langle f_{i_m}, f_{i_2} \rangle \\
\vdots & \vdots & \ddots & \vdots \\
\langle f_{i_1}, f_{i_m} \rangle & \langle f_{i_2}, f_{i_m} \rangle & \dots & \langle f_{i_m}, f_{i_m} \rangle \\
\end{array}
\right]
$$
\\

\noindent If we loose only one coefficient, then it can be seen that norm of error operator is given by norm of the corresponding frame vector. Therefore, in case of a uniform frame, the error is constant for each coefficient. 
\\

\noindent In case of losing two coefficients, $i^{th}$ and $j^{th}$, the error operator $ E_{i,j}$ is given by
$$ E_{i,j} (x)= x-\sum_{l \neq i,j} \langle x,f_{l} \rangle f_{l}= \langle x,f_{i} \rangle f_{i} + \langle x,f_{j} \rangle f_{j}$$

\noindent This condition is called {\em two-erasures}. A {\em two-uniform} frame is a frame which is uniform and $ \vert \vert E_{i,j} \vert \vert = \ constant \ \ \forall i \neq j $ . \\

\noindent We state the following characterization of two-uniform frames from \cite{2}.\\

\begin{thm}
Let $\lbrace f_{i} \rbrace_{i \in I}$ be a uniform $(n,k)$-frame. Then, it is two-uniform (equiangular) if and only if  $|\langle f_{i},f_{j} \rangle | = c_{n,k}$ for each $i,j \in I$ such that $i \neq j$ where 
$$ c_{n,k} = \sqrt{\dfrac{k(n-k)}{n^{2}(n-1)}}.$$\\
\end{thm}

\noindent A frame $\lbrace f_{i} \rbrace_{i \in I}$ such that $|\langle f_{i},f_{j} \rangle |= \ constant$  for each $i,j \in I$ for $i \neq j $ is called an {\em equiangular} frame.\\

\noindent In \cite{1}, it is established that when $\lbrace f_{i} \rbrace_{i \in I}$ is a uniform $(n,k)$-frame, then each vector $f_{i}$, $i\in I$ is of length $\sqrt{\frac{k}{n}}$. Also, it shows that two-uniform frames, when they exist, are optimal for two-erasures.\\

\noindent We now look at some basic concepts in Number Theory.\\

\noindent We know that $a \in \mathbb{Z}_n$ such that $gcd(a,n)=1$, is called a $quadratic$ $ residue$ of an odd prime $n$ if and only if $x^2 \equiv a \pmod n$ has a solution in $\mathbb{Z}_n$. Otherwise, $a$ is called a $quadratic$ $ nonresidue$ of $n$. Note that if $ a \equiv b \pmod{n}$, then $a$ is a quadratic residue (nonresidue) of $n$ if and only if $b$ is a quadratic residue (nonresidue) of $n$. Therefore, we  only look for residues in $\mathbb{Z}_{n}$. This is also called a $reduced$ $residue$ system. Note that the product of two quadratic residues or two quadratic nonresidues is a  quadratic residue in a reduced residue system of $n$.\\

\noindent Since $a^{n-1} \equiv 1 \pmod n$, then $a^{n-1}-1 \equiv (a^{\frac{n-1}{2}}-1)(a^{\frac{n-1}{2}}+1) \equiv 0 \pmod n$. \\

\noindent Thus, $a^{\frac{n-1}{2}} \equiv 1 \pmod n$ or $a^{\frac{n-1}{2}} \equiv -1 \pmod n$. \\

\noindent We can now state the following result also known as the {\em Euler's} criterion.\\

\begin{thm}

Let $n$ be an odd prime and $gcd(a,n)=1$. Then,

$$ a^{\frac{n-1}{2}}=(\dfrac{a}{n})_L = \left\{
\begin{array}{rcl}
\ 1 & & a\ is \ a \ quadratic \ residue\\
-1 & & a\ is \ a \ quadratic \ nonresidue\\
\end{array}\right.
$$
\\
\noindent where $()_L$ is called the {\em Legendre} symbol.
\end{thm}

\noindent It is easy to prove the following properties of Legendre symbol, \\

\noindent $1.\ (\frac{ab}{n})_L = (\frac{a}{n})_L (\frac{b}{n})_L .$\\

\noindent $2.\ (\frac{1}{n})_L =1 .$\\

\noindent Using modular arithmetics and the binomial theorem, we get 
$$(n-a)^{\frac{n-1}{2}} \equiv (-1)^{\frac{n-1}{2}}a^{\frac{n-1}{2}} \pmod{n}.$$

\noindent Let $n=2k+1$ be such that $n$ is prime and $k$ is odd. For such $k$ odd and $n$ prime, we equivalently say $n \equiv 3 \pmod{4}$. Let $a$ be in the reduced residue system of $n$. Thus by Euler's Criterion, $a$ is a quadratic residue (nonresidue) of $n$ if and only if $(n-a)$ is a quadratic nonresidue (residue) of $n$.

\begin{thm}
Let $n \equiv 3 \pmod 4$ be an odd prime. Then any reduced residue system $\pmod n$ contains $\frac{n-1}{2}$ quadratic residues and $\frac{n-1}{2}$ quadratic non-residues of $n$. One set of $\frac{n-1}{2}$ congruent quadratic residues is $ \lbrace 1^2, 2^2, \dots , \left(\frac{n-1}{2}\right)^2 \rbrace$.\\
\end{thm}

\noindent We now state the following result which provides a means for determining which primes have $2$ as a quadratic residue.\\

\begin{thm}
For an odd prime $n$, we have 

$$ (\dfrac{2}{n})_L = \left\{
\begin{array}{rcl}
\ 1 & & n \equiv \pm 1 \pmod{8}\\
-1 & & n \equiv \pm 3 \pmod{8}\\
\end{array}\right.
$$
\end{thm}

\noindent Let us now look at the case when $n \equiv 1 \pmod 4$ such that $n=4k^{2}+1$ is prime for odd $k$.

\begin{defn}

An element $a \in \mathbb{Z}_n$, such that $gcd(a,n)=1$, is called a quartic(biquadratic) residue of an odd prime $n$ if and only if $x^4 \equiv a \pmod n$ has a solution in $\mathbb{Z}_n$. Otherwise, $a$ is called a quartic(biquadratic) nonresidue of $n$.

\end{defn}

\noindent Note that every quartic residue is a quadratic residue. Also, product of two quartic residues is a quartic residue. By \cite{7}, we know that $a$ is a quartic residue of $n$ if $a^{\frac{n-1}{4}} \equiv 1 \pmod n$.\\

\noindent Let us denote the set of quadratic residues by $S_2$ and the set of quartic residues by $S_4$. Then it can be shown that for primes $n \equiv 1 \pmod 4$, $a \in S_2$ if and only if $(n-a) \in S_2$. Also, $a \in S_4$ if and only if $(n-a) \in S_2$. Therefore, $\overline{S}_4=S_{2}-S_{4}$ and $\lbrace S_4, \overline{S}_4 \rbrace$ forms a partion of $S_2$ such that $\arrowvert S_{4}\arrowvert = \arrowvert \overline{S}_4\arrowvert $.\\

\noindent Let $a \notin S_2$. Then it can be checked that $\lbrace aS_4,\overline{a}S_4 \rbrace $ forms a partition of  ${S_2}^{c}$ with $x \in aS_{4}$ if and only if $(n-x) \in \overline{a}S_{4}$.\\

\noindent Hence ${\mathbb{Z}_n}^{*}$ can be partioned into $\lbrace S_4, \overline{S_4},aS_4,\overline{a}S_4 \rbrace $ such that for every $x \in S_4$, we have $xS_4 = S_4$, $x \overline{S}_4 = \overline{S}_4$, $xaS_4 = aS_4 $ and $x \overline{a}S_4 = \overline{a}S_4 $.\\

\noindent Now let us look at the concept of difference sets.\\

\begin{defn}

\noindent A subset $H$ of a finite (additive) Abelian group $G$ is said to be a $(n,k,\lambda)$-difference set of $G$ if for some fixed natural number $\lambda$, every nonzero element of $G$ can be written as a difference of two elements of $H$ in exactly $\lambda$ ways, where $\arrowvert G \arrowvert = n$ and $\arrowvert H \arrowvert = k$.

\end{defn}

\noindent The followings tabulation on difference sets is discussed in \cite{8}. \\

{\em Type S}. (Singer difference sets). These are hyper planes in $PG(m,q)$, $q=p^r$. The parameteres are 
$$n=\frac{q^{m+1}-1}{q-1}, \ \ \ k=\frac{q^m -1}{q-1},\ \ \ \lambda=\frac{q^{m-1}-1}{q-1}. $$\\

{\em Type Q}. Let $n=p^r \equiv 3 \pmod 4 $. Then the quadratic residues of $n$ form a difference set with parameters
$$n=p^r=4t-1,\ \ \ k=2t-1, \ \ \ \lambda=t-1.$$

{\em Type $H_6$}. Let $n=4x^2+27$. There will exist a primitive root $r \pmod n$ such that $ Ind_r (3) \equiv 1 \pmod 6$. The residues $a_i \pmod n$ such that $Ind_r (a_i) \equiv 0,1 or 3 \pmod 6$ will form a difference set with 
$$n=4t-1,\ \ \ k=2t-1, \ \ \ \lambda=t-1.$$

{\em Type T} (Twin primes). Let $n$ and $n'=n+2$ be both primes. Then the collection of residues $\lbrace a_1, a_2, \dots, a_{m}, 0,n',2n',\dots, (n-1)n' \rbrace$ such that $(\frac{a_i}{p})_L =(\frac{a_i}{q})_L \ \ \forall i$ form a difference set $\pmod {nn'}$ with parameters
$$nn'=4t-1,\ \ \ k=2t-1, \ \ \ \lambda=t-1.$$  

\noindent Note that the types, {\em Q, $H_6$} and {\em T}, are {\em Hadamard} type difference sets.\\

{\em Type B}. Let $n=4x^2+1,\  x$ odd. Then the set of biquadratic (quartic) residues form a difference set with parameters
$$n=4x^2+1, \ \ \ k=x^2,\ \ \ \lambda = \frac{x^2-1}{4}. $$ 

{\em Type $B_0$}. Let $n=4x^2+9,\  x$ odd. Then the set of biquadratic (quartic) residues together with zero form a difference set with parameters
$$n=4x^2+9, \ \ \ k=x^2+3,\ \ \ \lambda = \frac{x^2+3}{4}. $$ 

{\em Type O}. Let $n=8x^2+1=64y^2+9,\ x$ and $y$ odd. Then the set of octic residues form a difference set with parameters
$$n=8x^2+1, \ \ \ k=x^2,\ \ \ \lambda = y^2. $$ 

{\em Type $O_0$}. Let $n=8x^2+49=64y^2+441,\  x$ odd and $y$ even. Then the set of octic residues together with zero form a difference set with parameters
$$n=8x^2+49, \ \ \ k=x^2+6,\ \ \ \lambda = y^2+7. $$ 

{\em Type $W_4$}. (Generalization of type T by Whiteman). Let $n$ and $n'=n+2$ be both primes such that $(n-1, n'-1)=4$. Define $d=(n-1)(n'-1)/4$. Let $g$ be a primitive root of both $n$ and $n'$. Then the collection of residues $\lbrace 1, g, g^2, \dots, g^{d-1}, 0,n',2n',\dots, (n-1)n'\rbrace$ form a difference set $\pmod {nn'}$ with parameters
$$nn',\ \ \ k=\frac{nn'-1}{4}, \ \ \ \lambda=\frac{nn'-5}{16}.$$  

\noindent More difference sets can be generated from a given difference sets. This can be seen from the following theorem.\\

\begin{thm}

\noindent A set $\lbrace n_1, n_2, \dots ,n_k \rbrace $ is a $(n,k,\lambda)$-difference set if and only if $\lbrace n_{1}+i, n_{2}+i, \dots ,n_{k}+i \rbrace $ for every $i \in \mathbb{Z}_n $ is a $(n,k,\lambda)$-difference set.\\

\end{thm}

\noindent Let $n$ be an odd prime and let $\lbrace n_1, n_2, \dots ,n_k \rbrace $ be a $(n,k,\lambda)$-difference set. Then $\lbrace \alpha n_1, \alpha n_2, \dots ,\alpha n_k \rbrace $ for every ${\mathbb{Z}_n}^{*}$ is also a $(n,k,\lambda)$-difference set as $\alpha$ is invertible and $n_{i}-n_{j} \equiv a_{ij} \Leftrightarrow \alpha n_{i}- \alpha n_{j} \equiv \alpha a_{ij}$. \\

\section{Cyclic Subspaces and Cyclic Frames}

\noindent Cyclic codes are one of the most useful codes in binary coding. 

\begin{defn}
A code $\mathcal{C} \in {\mathbb{Z}_2}^{n} $ is cyclic if $(x_{n-1}, x_{n-2}, \dots ,x_{1}, x_{0}) \in \mathcal{C}$ implies $(x_{n-2}, x_{n-3}, \dots ,x_{0}, x_{n-1}) \in \mathcal{C}$.
\end{defn} 

\noindent Thus $ \mathcal{C}$ is cyclic if and only if $\mathcal{C} \subseteq \dfrac{\mathcal{P}(x)}{\langle x^{n}-1 \rangle}$ is an ideal. These codes are efficient in detecting burst errors. A burst error of size $d$ is an $n$-tuple whose non-zero enteries are in a consecutive span of $d$ coordinates and no fewer.\\

\noindent Cyclic frames are inspired by the cyclic codes. We now look at the construction of cyclic equiangular frames.\\

\noindent Let $\lbrace e_{i} \rbrace_{i=1}^{n}$ be the standard orthonormal basis of $\mathbb{C}^{n}$. Let $S$ be the cyclic shift operator on $\mathbb{C}^{n}$ such that $Se_{i} = e_{i+1} \pmod n $ $ \forall \ i = 1,2,...,n-1$ and  $Se_{n} = e_{1}$. Then $S$ can be written as

$$S=\left[
\begin{array}{cccccc}
0 & 0 & 0 & \cdots & \cdots & 1 \\
1 & 0 & 0 & \cdots & \cdots & 0 \\
0 & 1 & 0 & \cdots & \cdots & 0 \\
\vdots & \vdots  & \ddots & \ddots &  & \vdots \\
\vdots & \vdots  & & \ddots & \ddots  & \vdots \\
0 & 0 & 0 & \cdots & 1 & 0 \\
\end{array}
\right]\in \mathbb{M}_n.
$$
\\

\noindent Note that $SS^{*}=S^{*}S=I$.\\

\begin{defn} 
A $k$-dimensional subspace $M$ of $\mathbb{C}^{n}$ is called cyclic if $M$ is shift-invariant, i.e. $S(M) \subseteq M $. \\
\end{defn}

\noindent Let $w=e^\frac{2\pi i}{n}$ be an $n^{th}$ root of unity. Define for each $i \in \lbrace 0,1,\ldots ,n-1 \rbrace$, $v_{i}\in \mathbb{C}^{n}$ as

$$v_{i}= \left[
\begin{array}{c}
1\\
\overline{w}^{i}\\ 
\overline{w}^{2i}\\
\vdots \\
\overline{w}^{(n-1)i}\\
\end{array}
\right]
$$
\\
\noindent Note that 

\begin{equation}
Sv_{i}= w^{-i}v_{i}= \overline{w}^{i}v_{i}. 
\end{equation}

\noindent It can be shown that $\langle v_{i},v_{j}\rangle=0$ $\forall$ $i\neq j$, and $\parallel v_{i} \parallel = \sqrt{n}$. Thus $\lbrace\dfrac{1}{\sqrt{n}} v_i \rbrace_{i=0}^{n-1}$ is an orthonormal basis for $\mathbb{C}^{n}$.\\

\noindent Let $\lbrace f_0, f_1, \dots , f_{n-1} \rbrace $ be $(n,k)$-frame and define

$$V= \left[
\begin{array}{c}
{f_{0}}^{*}\\ 
{f_{1}}^{*}\\
\vdots \\
{f_{n-1}}^{*}\\
\end{array}
\right].
$$
\\

\noindent Then $V$ is an isometry.

\begin{defn}
The frame $\lbrace f_0, f_1, \dots , f_{n-1} \rbrace $ is called cyclic if and only if range$(V)$ is shift-invariant.\\
\end{defn}

\noindent Now considering the cyclic shift operator $S$ and $k$-dimensional subspace $M$, we know $S(M)\subseteq M$. By above proposition, there exists $I \subseteq \lbrace 1,2,\dots, n \rbrace $ such that $M= span\lbrace v_{i}:\  i\in I\rbrace$. \\

\noindent By using (1), we can summarize this in the following theorem.\\

\begin{thm}
Let $M$ be a subspace of $\mathbb{C}^{n}$. Then $M$ is $S$-invariant if and only if $\exists \  I\subseteq \lbrace 0,1,2,\ldots,n-1\rbrace$ such that $M= span\lbrace v_{i}:\  i\in I\rbrace$.\\
\end{thm}

\noindent Let $\lbrace n_{1}, n_{2},\ldots,n_{k}\rbrace \subseteq \lbrace0,1,2,\ldots,n-1\rbrace$ and $f_j \in \mathbb{C}^{k}$ such that for each $ j \in \lbrace 0,1,\dots,n-1 \rbrace$,

$$f_{j}= \dfrac{1}{\sqrt{n}}\left[
\begin{array}{c}
w^{jn_{1}}\\ 
w^{jn_{2}}\\
\vdots \\
w^{jn_{k}}\\
\end{array}
\right].
$$
\\
\noindent \noindent where $w$ is a primitive $n^{th}$ root of unity. These vectors form a $(n,k)$-frame. In \cite{3}, these are called the {\em harmonic frames}. For this frame,\\

$$V=\dfrac{1}{\sqrt{n}}\left[
\begin{array}{cccccc}
1 & 1 & 1 & \cdots & \cdots & 1 \\
\overline{w}^{n_{1}} & \overline{w}^{n_{2}} & \overline{w}^{n_{3}} & \cdots & \cdots & \overline{w}^{n_{k}} \\
\overline{w}^{2n_{1}} & \overline{w}^{2n_{2}} & \overline{w}^{2n_{3}} & \cdots & \cdots & \overline{w}^{2n_{k}} \\
\vdots & \vdots  & \vdots & \ddots &  & \vdots \\
\vdots & \vdots  & \vdots &  & \ddots  & \vdots \\
\overline{w}^{(n-1)n_{1}} & \overline{w}^{(n-1)n_{2}} & \overline{w}^{(n-1)n_{3}} & \cdots & \cdots & \overline{w}^{(n-1)n_{k}} \\
\end{array}
\right]
$$
\\
\noindent and $V:\mathbb{C}^k \rightarrow \mathbb{C}^n $ is an isometry. Let $M$ be the range of $V$.\\

\noindent Note that for any choice of $\lbrace n_1, n_2, \ldots , n_k \rbrace \subseteq \lbrace 0,1,\ldots , n-1 \rbrace $, $V$ is {\em shift-invariant}. Hence every harmonic frame is {\em cyclic}.\\

\noindent Note that $\lbrace f_{0},f_{1},\ldots,f_{n-1}\rbrace$ is a $\textsl{Parseval}$ frame for $\mathbb{C}^{k}$. This family of frames was introduced in \cite{3}, but we should also notice their cyclic nature. \\

\noindent Let $M$ be an $S$-invariant subspace of $\mathbb{C}^{n}$ and let $P_{M} : \mathbb{C}^{n}\rightarrow M$ be the orthogonal projection. We now prove the following theorem. \\

\begin{thm}
Let $M$ be a subspace of $\mathbb{C}^{n}$ such that $S(M)\subseteq M$ with orthogonal projection $P_{M}$. Then $SP_{M}=P_{M}S$.
\end{thm}

\begin{proof}
Let $v \in \mathbb{C}^{n}$ such that 
$$v=\sum_{v_{i}\in M} \alpha_{i} v_{i} + \sum_{v_{i}\notin M} \alpha_{i} v_{i}.$$
\noindent Then
$$ 
\begin{array}{lcl}
SP_{M}(v)&=& S(\sum _{v_{i}\in M} \alpha_{i} v_{i})\\
 &=&\sum _{v_{i}\in M} \alpha_{i} \overline{w}^{i} v_{i}\\
 &=& P_{M}(\sum _{v_{i}\in M} \alpha_{i} \overline{w}^{i} v_{i}+\sum _{v_{i}\notin M} \alpha_{i} v_{i})\\
 &=& P_{M}S(v).
\end{array}$$
\end{proof}

\noindent Consider $A=(a_{i,j})_{i,j}$ be in the commutant of $S$. Then 
$$
\begin{array}{lcl}
a_{i,j} & = & \langle Ae_{j}, e_{i}\rangle \\
 & = & \langle ASe_{j-1}, Se_{i-1}\rangle \\
 & = & \langle SAe_{j-1}, Se_{i-1}\rangle \\
 & = & \langle S^{*}SAe_{j-1}, e_{i-1}\rangle \\
 & = & \langle Ae_{j-1}, e_{i-1}\rangle \\
 & = & a_{i-1,j-1}
\end{array}$$

\noindent Therefore, every $A$ such that $AS=SA$ is of the form 

$$A=a_{0}I+a_{1}S+a_{2}S^{2}+ \ldots + a_{n-1}S^{n-1}$$ \\

\noindent for some constants $a_{0},a_{1},a_{2},\ldots,a_{n-1}$. So $P_{M}$ can be written as 

$$P_{M}=\left[
\begin{array}{cccccc}
a_{0} & a_{n-1} & a_{n-2} & \cdots & \cdots & a_{1} \\
a_{1} & a_{0} & a_{n-1} & \cdots & \cdots & a_{2} \\
a_{2} & a_{1} & a_{0} & \cdots & \cdots & a_{3} \\
\vdots & \vdots  & \ddots & \ddots &  & \vdots \\
\vdots & \vdots  & & \ddots & \ddots  & \vdots \\
a_{n-1} & a_{n-2} & a_{n-3} & \cdots & \ddots & a_{0} \\
\end{array}
\right]
$$
\\
where $P_{M}=VV^{*}=\left( \langle f_{j},f_{i}\rangle \right)_{i,j} $ is {\em Toeplitz}. Moreover, $P_{M}$  is {\em circulant} due to its cyclic nature. Note that this frame forms an ordered collection as changing the order of $f_{i}$'s disturbs the Toeplitz structure of $P_{M}$.\\

\begin{thm}
Let $\lbrace h_1,h_2,\dots , h_n \rbrace $ be a cyclic $(n,k)$-frame. Then $\exists$ a unitary $U$ and $\lbrace n_1,n_2,\ldots,n_k \rbrace \subseteq \lbrace 0,1,\dots , n-1 \rbrace $ such that $ Uh_i = f_i \ \forall i$, where  

$$f_{j}= \dfrac{1}{\sqrt{n}}\left[
\begin{array}{c}
w^{jn_{1}}\\ 
w^{jn_{2}}\\
\vdots \\
w^{jn_{k}}\\
\end{array}
\right].
$$ 

\noindent defines a harmonic $(n,k)$-frame.\\
\end{thm}

\begin{proof}
\noindent Let $\lbrace h_1,h_2,\dots , h_{n} \rbrace $ be a cyclic $(n,k)$-frame. Then $P_H =\left(\langle h_j,h_i \rangle \right)_{i,j}$ is Toeplitz and can be expressed as a polynomial $p(S)$. We know that the eigenvalues of $S$ are $\lbrace w^l:l=0,1,\dots,n-1 \rbrace $ where $w$ is the primitive $n^{th}$ root of unity. Thus, the eigenvalues of $P_H$ are given by  $\lbrace p(w^l):l=0,1,\dots,n-1 \rbrace $. Since $P_H$ is also a projection and trace$(P_H)=k$, therefore there are exactly $k$ $1$'s and $(n-k)$ $0$'s. Let $p(w^l)=1$ for $l=n_1,n_2,\dots,n_k$ and $p(w^l)=0$ otherwise.\\

\noindent Since the eigenvectors of $S$ are also the eigenvectors of $P_H = p(S)$, therefore for each $j=1,2,\dots, n$ the eigenvector of $P_H$ are given by 

$$v_{j}= \left[
\begin{array}{c}
1\\
w^{j(n-1)}\\ 
w^{j(n-2)}\\
\vdots \\
w^{j}\\
\end{array}
\right].
$$ 

\noindent For each $j=1,2,\dots, n$, define

$$f_{j}= \dfrac{1}{\sqrt{n}}\left[
\begin{array}{c}
w^{jn_{1}}\\ 
w^{jn_{2}}\\
\vdots \\
w^{jn_{k}}\\
\end{array}
\right].
$$ 

\noindent Then $\lbrace f_1,f_2,\dots, f_n \rbrace $ is a cyclic $(n,k)$-frame. Let $P_F=\left(\langle f_j,f_i \rangle \right)_{i,j}$. For some $r$, consider the $i^{th}$ entry of the vector $P_F(v_r)$, $i.e$,

$$ \sum_{j=1}^{n} \langle f_j, f_i \rangle v_{rj} =  \sum_{j=1}^{n} \sum_{t=1}^{k} w^{(j-i)n_t} w^{r(-j+1)} =\sum_{t=1}^{k} w^{(r-in_t)} \sum_{j=1}^{n} w^{(n_t -r)j}.$$

\noindent Thus the vector $v_r$ is a zero vector exactly when $r \in \lbrace n_1,n_2,\dots,n_k \rbrace$. Therefore, both $P_F$ and $P_M$ have exactly same eigenvalues and eigenvectors. Hence $P_H=P_F$. Thus by \cite{2}, the cyclic frame $\lbrace h_1,h_2,\dots , h_{n} \rbrace $ is unitarily equivalent to the frame  $\lbrace f_1,f_2,\dots, f_n \rbrace $.\\

\end{proof}

\noindent Hence it suffices to only consider the cyclic frames of type $\lbrace f_1, f_2, \dots , f_n \rbrace $ defined above. 
\\

\noindent We will now try to find optimal cyclic frames for two-erasures indexed by the subset $I \subseteq \lbrace 0,1,2,\ldots, n-1 \rbrace$, where order of $I$ is $k$. \\
\\
\noindent Since for optimal equiangular cyclic frames we have $|\langle f_{j},f_{i} \rangle |=constant$ and $P_{M}$ is a circulant matrix, therefore, the problem to find optimal equiangular cyclic frames is now reduced to showing that 
$$|\langle f_{i},f_{0} \rangle | = |\langle f_{j},f_{0} \rangle | \ \ \ \ \forall \  i \neq j,\  i \neq 0 ,\  j \neq 0 .$$

\noindent Hence, we need to find a subset $I \subseteq \lbrace 0,1,2,\ldots, n-1 \rbrace$ such that for each  $j \in \lbrace 2,3,\ldots, n-1 \rbrace$, the \textit{absolute condition} is satisfed, \textit{i.e.},
$$ |\sum_{i=1}^{k} w^{n_{i}}| = |\sum_{i=1}^{k} w^{jn_{i}}|.$$
\\
\noindent In the rest of this paper, we will establish some conditions on $n$ and $k$ to show the existence of equiangular cyclic $(n,k)$-frames defined as above.\\

\section{Equiangular cyclic frames}

\noindent We now study the possible selections of the set $I \subseteq \mathbb{Z}_n$ in order to generate equiangular cyclic $(n,k)$-frames. We note that in order to determine an equiangular cyclic $(n,k)$-frame as above, we need to determine the frame vector $f_{1}$ only. Therefore, we will call the vector $f_{1}$ as the $generator$ of the frame. \\

\noindent Since the absolute value of a sum does not change when the entries are permuted, thus we will mainly consider the set $f'_{j}=\dfrac{1}{\sqrt{n}}\lbrace w^{jn_{1}},w^{jn_{2}},\ldots,w^{jn_{k}}\rbrace$ of the entries in vector 

$$f_{j}=\dfrac{1}{\sqrt{n}}\left[
\begin{array}{c}
w^{jn_{1}}\\ 
w^{jn_{2}}\\
\vdots \\
w^{jn_{k}}\\
\end{array}
\right].
$$

\noindent Therefore although we will state the results for the vectors $\lbrace f_j \rbrace_{j \in J} $, but it would suffice to prove it for the sets  $\lbrace f'_j \rbrace_{j \in J} $.\\

\noindent These frames were studied in \cite{9} where the following result is proved for complex MWBE(Maximum Welch Bound Equality) codebooks. The construction of these codebooks show that they are the same as equiangular cyclic frames and therefore, we state the theorem for equiangular cyclic frame and provide a slightly different proof.\\

\begin{thm}
\noindent The collection $\lbrace f_{0},f_{1},f_{2},\ldots,f_{n-1}\rbrace $ is an equiangular cyclic $(n,k)$-frame if and only if the set $\lbrace n_1,n_2,\dots,n_k\rbrace$ is a $(n,k,\lambda)$-difference set, where
$$f_{j}=\dfrac{1}{\sqrt{n}}\left[
\begin{array}{c}
w^{jn_{1}}\\ 
w^{jn_{2}}\\
\vdots \\
w^{jn_{k}}\\
\end{array}
\right].
$$
\end{thm}

\begin{proof}
Let the collection $\lbrace f_{0},f_{1},f_{2},\ldots,f_{n-1}\rbrace $ be an equiangular cyclic $(n,k)$-frame. Then by \cite{2}, we know that for every $l \neq 0$, we have

$$\arrowvert \langle f_0, f_l \rangle \arrowvert^2 =\dfrac{1}{n^2 }\sum_{i,j=1}^{k} w^{l(n_i - n_j)}=c^2$$

\noindent where $c= \sqrt{\frac{k(n-k)}{n^2 (n-1)}}$. \\

\noindent Let $a_r$ be the order of the set $\lbrace (i,j)\ |\ n_i-n_j \equiv r \pmod n \rbrace$. Then $a_0=k$. Clearly, 

$$\dfrac{1}{n^2}\sum_{i,j=1}^{k} w^{l(n_i - n_j)}=\dfrac{1}{n^2}\sum_{t=0}^{n-1}a_t w^{lt}=c^2. $$ \\ 

\noindent Let $p(z)=a_0+a_1 z+\dots+a_{n-1} z^{n-1}$. Then $p(w^l)=n^2 c^2 \ \ \forall \ l \neq 0$ and $p(1)=a_0+a_1+\dots+a_{n-1}$.\\

\noindent Consider the $n \times n $ matrix $U=(w^{ij})_{i,j}$. Then  $U^{*}=(w^{-ij})_{i,j}$ and 

$$U^{*}U=\left(\sum_{t=0}^{n-1}w^{-it}w^{tj}\right)_{i,j}=\left(\sum_{t=0}^{n-1}w^{(j-i)t}\right)_{i,j}=nI.$$\\

\noindent Therefore, 
$$
U \left[
\begin{array}{c}
a_0\\
a_1\\
a_2\\
\vdots\\
a_{n-1}\\
\end{array}
\right]=\left[
\begin{array}{c}
a_0+a_1+\dots+a_{n-1}\\
p(w)\\
p(w^2)\\
\vdots\\
p(w^{n-1)}\\
\end{array}
\right]=\left[
\begin{array}{c}
a_0+a_1+a_2+\dots+a_{n-1}\\
n^2 c^2\\
n^2 c^2\\
\vdots\\
n^2 c^2\\
\end{array}
\right]
$$

\noindent So, we get, 
$$
n\left[
\begin{array}{c}
a_0\\
a_1\\
a_2\\
\vdots\\
a_{n-1}\\
\end{array}
\right]=U^{*}U\left[
\begin{array}{c}
a_0\\
a_1\\
a_2\\
\vdots\\
a_{n-1}\\
\end{array}
\right]=U^{*}\left[
\begin{array}{c}
a_0+a_1+a_2+\dots+a_{n-1}\\
n^2 c^2\\
n^2 c^2\\
\vdots\\
n^2 c^2\\
\end{array}
\right].
$$

\noindent Thus,
$$
n\left[
\begin{array}{c}
a_0\\
a_1\\
a_2\\
\vdots\\
a_{n-1}\\
\end{array}
\right]=\left[
\begin{array}{c}
a_0+a_1+a_2+\dots+a_{n-1}+(n-1)n^2 c^2\\
a_0+a_1+a_2+\dots+a_{n-1}-n^2 c^2 \\
a_0+a_1+a_2+\dots+a_{n-1}-n^2 c^2 \\
\vdots \\
a_0+a_1+a_2+\dots+a_{n-1}-n^2 c^2\\
\end{array}\right].
$$

\noindent Hence, $na_r = a_0+a_1+a_2+\dots+a_{n-1}-n^2 c^2$, which is independent of $r$, therefore, 
$$a_1 = a_2 = \dots =a_{n-1} = \lambda $$

\noindent for some constant $\lambda$. So, $n \lambda=k+(n-1)\lambda - n^2 c^2$. Thus, 

$$\lambda = k-n^2 c^2 = k-\dfrac{k(n-k)}{n-1}=\dfrac{k(k-1)}{n-1}.$$  

\noindent Hence, $\lbrace n_1,n_2,\dots,n_k\rbrace$ forms a difference set.\\

\noindent Conversely, let $\lbrace n_1,n_2,\dots,n_k\rbrace$ be a difference set. For each non-zero $j$, define 
$$f_{j}=\dfrac{1}{\sqrt{n}}\left[
\begin{array}{c}
w^{jn_{1}}\\ 
w^{jn_{2}}\\
\vdots \\
w^{jn_{k}}\\
\end{array}
\right].
$$

\noindent Then \cite{2} showed that the collection $\lbrace f_{0},f_{1},f_{2},\ldots,f_{n-1}\rbrace $ is equiangular if and only if 

$$\arrowvert \langle f_i , f_j \rangle \arrowvert =c$$

\noindent for all $i \neq j$. Consider

$$
\begin{array}{lcl}
\arrowvert \langle f_i , f_j \rangle \arrowvert^2 & = & \arrowvert \langle f_i , f_j \rangle \arrowvert\overline{\arrowvert \langle f_i , f_j \rangle \arrowvert}\\
 & = & \frac{1}{n^2} \left(\sum_{l=1}^{k}w^{in_l -jn_l} \right) \left( \sum_{m=1}^{k}w^{jn_m -in_m}\right)\\
 & = & \frac{1}{n^2} \left( \sum_{l,m=1}^{k}w^{(i-j)(n_l - n_m)} \right)\\
 & = & \frac{1}{n^2} \left( k+ \sum_{l \neq m}w^{(i-j)(n_l - n_m)} \right)\\
 & = & \frac{1}{n^2} \left( k+ \sum_{r=1}^{n-1} a_r w^{(i-j)r} \right)\\
 & = & \frac{1}{n^2} \left( k+ \sum_{r=1}^{n-1} \lambda w^{(i-j)r} \right)\\
 & = & \frac{1}{n^2} \left( k+ \lambda \sum_{r=1}^{n-1} w^{(i-j)r} \right)\\
 & = & \frac{1}{n^2} \left( k+ \lambda (-1) \right)\\
 & = & \frac{1}{n^2} \left( k- \lambda \right)\\
 & = & c^2.\\
\end{array}
$$

\noindent Hence, as the absolute condition is satisfied, therefore the collection $\lbrace f_{i}\rbrace_{i=0}^{n-1} $ is an equiangular cyclic $(n,k)$-frame. \\

\end{proof}

\noindent Since $\lambda = \frac{k(k-1)}{n-1}$ must be an integer, therefore we get the following corollary which gives a necessary condition for the existence of a  $(n,k,\lambda)$-difference set.\\

\begin{cor}
Let there exist a $(n,k,\lambda)$-difference set, then $n-1$ must divide $k(k-1)$.
\end{cor}

\noindent We now look at the following results which are obtained independently of \cite{9}. These examine equiangular cyclic frames without involving difference sets. A close observation of the results from \cite{9} with the following results reveal some very interesting properties of difference sets which might not be that obvious by definition.\\

\noindent We start with an example of an equiangular cyclic $(7,3)$-frame and depict the use of the {\em absolute condition}.\\

\noindent For $n=7$ and $k=3$, let us choose 

$$f_{1}=\dfrac{1}{\sqrt{7}}\left[
\begin{array}{c}
w\\ 
w^{2}\\
w^{3}\\
\end{array}
\right].
$$

\noindent So $f'_{1}=\{ w,w^{2},w^{3}\}$. We can check that the  absolute condition fails as 

$$|\langle f_1, f_0 \rangle|=|w+w^{2}+w^{3}| \neq |w^{2}+w^{4}+w^{6}|=|\langle f_2, f_0 \rangle|.$$

\noindent Thus, the above chosen $f_{1}$ does not generate an equiangular cyclic $(7,3)$-frame.\\ 

\noindent However, let us now choose 

$$f_{j}=\dfrac{1}{\sqrt{7}}\left[
\begin{array}{c}
w\\ 
w^{2}\\
w^{4}\\
\end{array}
\right].
$$\\

\noindent So $f'_{1}=\{ w,w^{2},w^{4}\}$.  Then we can check that $f'_{i}=\lbrace w,w^{2},w^{4}\rbrace$ for $i=1,2,4$ and $f'_{i}=\lbrace w^{3},w^{5},w^{6}\rbrace$ for $i=3,5,6$, and

$$|w+w^{2}+w^{4}| = |w^{3}+w^{5}+w^{6}|$$\\

\noindent Thus the absolute condition is satisfied and hence $f_{1}$ generates such an equiangular cyclic $(7,3)$-frame. \\

\noindent We  look at the first theorem which demonstrates the relation between $n,k$ and the roots of unity required to generate an equiangular cyclic $(n,k)$-frame.\\

\begin{thm}

\noindent Let $ \lbrace n_{1},n_{2},\ldots,n_{k} \rbrace \subseteq \mathbb{Z}_{n}$ such that $\forall \  j \in \lbrace 2,3,\ldots,n-1 \rbrace$, $\exists \  l_{j}\in \mathbb{Z}_{n}$ and a permutation $\pi_{j}$ of $\mathbb{Z}_{k}$ such that either
$$(i)\ \  jn_{i}-n_{\pi_{j}(i)}\equiv l_{j}\pmod{n}\ \ \  \forall i=1,2,\ldots,k $$
or,
$$(ii)\ \ jn_{i}+n_{\pi_{j}(i)}\equiv l_{j}\pmod{n}\ \ \  \forall i=1,2,\ldots,k  $$
\\
\noindent Then the collection $\lbrace f_{0},f_{1},f_{2},\ldots,f_{n-1}\rbrace $ is an equiangular cyclic $(n,k)$-frame with 

$$f_{j}=\dfrac{1}{\sqrt{n}}\left[
\begin{array}{c}
w^{jn_{1}}\\ 
w^{jn_{2}}\\
\vdots \\
w^{jn_{k}}\\
\end{array}
\right].
$$

\end{thm}

\begin{proof}
 Let  $ \lbrace n_{1},n_{2},\ldots,n_{k} \rbrace \subseteq \mathbb{Z}_{n}$ be chosen as above. Then  $\lbrace f_{0},f_{1},f_{2},\ldots,f_{n-1}\rbrace $ will be an equiangular cyclic $(n,k)$-frame if and only if the absolute condition is satisfied, $\textit{i.e.}$ for each $j\in \lbrace 2,3,\ldots,n-1\rbrace$
 
$$ |\sum_{i=1}^{k} w^{n_{i}}| = |\sum_{i=1}^{k} w^{jn_{i}}|.$$\\

\noindent Firstly, choose $j$ such that $(i)$ is satisfied. Then for  $j\in \lbrace 2,3,\ldots,n-1\rbrace$, 

$$f_{j}=\dfrac{1}{\sqrt{n}}\left[
\begin{array}{c}
w^{jn_{1}}\\ 
w^{jn_{2}}\\
\vdots \\
w^{jn_{k}}\\
\end{array}
\right].
$$

\noindent Then $f'_{j}=\{ w^{jn_{1}}, w^{jn_{2}},\ldots, w^{jn_{n-1}} \}.$ \\

\noindent Then,
$$w^{jn_{m}}=w^{l_{j}+n_{\pi_{j}}(m)}=w^{l_{j}}w^{n_{\pi_{j}}(m)}.$$ 

\noindent Therefore,

$$\begin{array}{l}
 |\sum_{i=1}^{k} w^{jn_{i}}|\\
= |\sum_{i=1}^{k} w^{l_{j}}w^{n_{\pi_{j}(i)}}|\\
= | w^{l_{j}}\sum_{i=1}^{k} w^{n_{\pi_{j}(i)}}|\\
= | w^{l_{j}}||\sum_{m=1}^{k} w^{n_{m}}|\\
= |\sum_{m=1}^{k} w^{n_{m}}|\\
\end{array}
$$\\

\noindent Now, choose $j$ such that $(ii)$ is satisfied. Then,

$$w^{jn_{m}}=w^{l_{j}-n_{\pi_{j}}(m)}=w^{l_{j}}\overline{w^{n_{\pi_{j}}(m)}}.$$ 

\noindent Therefore,

$$
\begin{array}{l}
 |\sum_{i=1}^{k} w^{jn_{i}}|\\
= |\sum_{i=1}^{k} w^{l_{j}}\overline{w^{n_{\pi_{j}(i)}}}|\\
= | w^{l_{j}}\sum_{i=1}^{k} \overline{w^{n_{\pi_{j}(i)}}}|\\
= | w^{l_{j}}||\overline{\sum_{m=1}^{k} w^{n_{m}}}|\\
= |\sum_{m=1}^{k} w^{n_{m}}|\\
\end{array}
$$
\\
\noindent Thus the absolute condition is satisfied for all $j$. Hence, $\lbrace f_{0},f_{1},f_{2},\ldots,f_{n-1}\rbrace $ is an equiangular cyclic $(n,k)$-frame.
\end{proof}

\noindent As shown in \cite{9}, a set $ \lbrace n_{1},n_{2},\ldots,n_{k} \rbrace \subseteq \mathbb{Z}_{n}$ generates an equiangular cyclic $(n,k)$-frame if and only if the set is a difference set, therefore we have the following corollary.\\

\begin{cor}

Any set $ \lbrace n_{1},n_{2},\ldots,n_{k} \rbrace \subseteq \mathbb{Z}_{n}$ satisfying the properties of Theorem $4.2$ must be a difference set.\\

\end{cor}

\noindent More frames can be developed from a given frame. We use the fact that the sum of all roots of unity is zero to prove the following theorem.\\

\begin{thm}
Let $w$ be a primitive $n^{th}$ root of unity. Let $n_{1},n_{2},\ldots,n_{k} \in \mathbb{Z}_{n}$ such that 

$$f_{1}=\dfrac{1}{\sqrt{n}}\left[
\begin{array}{c}
w^{n_1}\\ 
w^{n_2}\\
\vdots \\
w^{n_k}\\
\end{array}
\right]
$$ \\

\noindent generates an equiangular cyclic $(n,k)$-frame. Then the remaining roots of unity generate an equiangular cyclic $(n,n-k)$-frame generated by
 
$$h_{1}=\dfrac{1}{\sqrt{n}}\left[
\begin{array}{c}
w^{n_{k+1}}\\ 
w^{n_{k+2}}\\
\vdots \\
w^{n_{n}}\\
\end{array}
\right]
$$\\
\end{thm}

\begin{proof}
Let  $ \lbrace n_{k+1},n_{k+2},\ldots,n_{n} \rbrace \subseteq \mathbb{Z}_n$ such that $\{ w^{n_{k+1}}, w^{n_{k+2}},\ldots, w^{n_{n}} \}$ forms the set of remaining $n^{th}$ roots of unity. Since 

$$w^{n_{1}}+w^{n_{2}}+\ldots+w^{n_{n}}=0,$$

\noindent therefore,

$$w^{jn_{1}}+w^{jn_{2}}+\ldots+w^{jn_{n}}=0.$$

\noindent As $n$ is prime, then

$$w^{jn_{r}}=w^{jn_{s}} \iff n_{r} \equiv n_{s} \pmod{n}.$$

\noindent So,

$$w^{jn_{1}}+w^{jn_{2}}+\ldots+w^{jn_{k}}=-w^{jn_{k+1}}-w^{jn_{k+2}}-\ldots+w^{jn_{n}}\ \ \forall \  1\leq j\leq n-1.$$\\

\noindent Then,

$$|w^{n_{1}}+w^{n_{2}}+\ldots+w^{n_{k}}|=|w^{n_{k+1}}+w^{n_{k+2}}+\ldots+w^{n_{n}}| $$

\noindent and,

$$|w^{jn_{1}}+w^{jn_{2}}+\ldots+w^{jn_{k}}|=|w^{jn_{k+1}}+w^{jn_{k+2}}+\ldots+w^{jn_{n}}|.$$\\

\noindent But $\lbrace f_{i}\rbrace$ is an equiangular cyclic $(n,k)$-frame, therefore,

$$|w^{n_{1}}+w^{n_{2}}+\ldots+w^{n_{k}}|=|w^{jn_{1}}+w^{jn_{2}}+\ldots+w^{jn_{k}}|.$$\\

\noindent Therefore the absolute condition is satisfied, i.e., 

$$|w^{n_{k+1}}+w^{n_{k+2}}+\ldots+w^{n_{n}}|=|w^{jn_{k+1}}+w^{jn_{k+2}}+\ldots+w^{jn_{n}}| .$$\\

\noindent So $\lbrace g_{0},g_{1},\ldots,g_{n-1}\rbrace$ generated by $g_{1}=( w^{n_{k+1}}, w^{n_{k+2}},\ldots, w^{n_{n}})$ is an equiangular cyclic $(n,k+1)$-frame.

\end{proof}

\noindent Then by Theorem $3.4$ we get that the frame generated by $g'_{1}=( 1,w^{3},w^{5},w^{6} )$ is an equiangular cyclic $(7,4)$-frame.\\

\begin{cor}
Let a given set $\lbrace n_{1},n_{2},\ldots,n_{k} \rbrace \subseteq \mathbb{Z}_{n}$ be a $(n,k,\lambda)$-difference set. Then the complement must also be a $(n,n-k,\bar{\lambda})$-difference set, where 

$$\bar{\lambda}=\dfrac{(n-k)(n-k-1)}{n-1}. $$

\end{cor}

\noindent Since $\bar{\lambda} = \frac{(n-k)(n-k-1)}{n-1}$ must be an integer, therefore we get the following corollary which now gives a stronger necessary condition for the existence of a  $(n,k,\lambda)$-difference set.\\

\begin{cor}
Let there exist a $(n,k,\lambda)$-difference set, then $n-1$ must divide the two quanties $k(k-1)$ and $(n-k)(n-k-1)$.\\
\end{cor}

\noindent Note that it is true for any equiangular cyclic $(n,k)$-frame which is generated by $n^{th}$ roots of unity as above. It can also be deduced from above that we always have an equiangular cyclic $(n,1)$-frame, and hence, an equiangular cyclic $(n,n-1)$-frame.\\

\noindent Look at the example of an equiangular cyclic $(7,3)$-frame. It should be noticed that the set $\lbrace 1,2,4 \rbrace$ is the set of quadratic residues of $7$. Also, these are the powers of $w$ which generate an equiangular cyclic $(7,3)$-frame. \\

\noindent We  now generalize the above observation from the case of an equiangular cyclic $(7,3)$-frame.\\

\begin{thm}
Let $n$ be a prime integer such that $n=2k+1$, where $k$ is odd. For each $j \in \lbrace 0,1,\ldots,n-1 \rbrace$, define 

$$f_{j}=\dfrac{1}{\sqrt{n}}\left[
\begin{array}{c}
w^{j1^2}\\ 
w^{j2^2}\\
\vdots \\
w^{jk^2}\\
\end{array}
\right]
$$\\

\noindent Then the collection $\lbrace f_{0}, f_{1},\ldots,f_{n-1} \rbrace$ is an equiangular cyclic $(n,k)$-frame.\\
\end{thm}

\begin{proof}
Let $m$ be in the reduced residue system of $n$. Since in a reduced residue system, the product of two quadratic residues (non-residues) is a quadratic residue and product of a quadratic non-residue with a quadratic residue is a quadratic non-residue, thus $jm^{2}$ is a quadratic residue (nonresidue) if and only if $j$ is a quadratic residue (nonresidue). Then $ \lbrace -{1^2}, -{2^2}, \ldots, -{k^2} \rbrace$ be the set of quadratic nonresidues of  $n$. Since the set of all quadratic residues is closed with respect to multiplication, therefore, $f'_{j} = \lbrace w^{1^{2}},w^{2^{2}}, \ldots ,w^{k^{2}}\rbrace $ when j is a quadratic residue, and, $ f'_{j}=\lbrace w^{-{1^2}},w^{-{2^2}}, \ldots ,w^{-{k^2}}\rbrace$ when $j$ is a quadratic nonresidue. Since $a$ is a quadratic residue (nonresidue) of $n$ if and only if $(n-a)$ is a quadratic nonresidue (residue) of $n$, therefore \\
$$
\begin{array}{l}
|w^{1^{2}}+w^{2^{2}}+\ldots+w^{k^{2}}|\\
=|\overline{w^{1^{2}}+w^{2^{2}}+\ldots+w^{k^{2}}}|\\
=|\overline{w^{1^{2}}}+\overline{w^{2^{2}}}+\ldots+\overline{w^{k^{2}}}|\\
=|\overline{w}^{1^{2}}+\overline{w}^{2^{2}}+\ldots+\overline{w}^{k^{2}}|\\
=|w^{n-1^{2}}+w^{n-2^{2}}+\ldots+w^{n-k^{2}}|\\
=|w^{-{1^2}}+w^{-{2^2}}+\ldots+w^{-{k^2}}|\\
\end{array}
$$

\noindent Since the absolute condition is satisfied, therefore $\lbrace f_{0}, f_{1},\ldots,f_{n-1} \rbrace$ is an equiangular cyclic $(n,k)$-frame.\\
\end{proof}

\noindent We now look at another way of generating frames from a given equiangular cyclic $(n,k)$-frame.\\

\begin{cor}
Let $n$ be a prime integer such that $n=2k+1$, where $k$ is odd. Let

$$f_{j}=\dfrac{1}{\sqrt{n}}\left[
\begin{array}{c}
w^{j1^2}\\ 
w^{j2^2}\\
\vdots \\
w^{jk^2}\\
\end{array}
\right]
$$\\

\noindent for each $j \in \lbrace 0,1,\ldots,n-1 \rbrace$ be the equiangular cyclic $(n,k)$-frame constructed as above.  Then 

$$Re(w^{1^{2}}+w^{2^{2}}+ \ldots +w^{k^{2}})=-\dfrac{1}{2}.$$

\noindent Consequently, 

$$g_{1}=\dfrac{1}{\sqrt{n}}\left[
\begin{array}{c}
1\\
w^{1^2}\\ 
w^{2^2}\\
\vdots \\
w^{k^2}\\
\end{array}
\right]
$$

\noindent generates an equiangular cyclic $(n,k+1)$-frame.\\
\end{cor}

\begin{proof}
Since all $n^{th}$ roots of unity sum up to zero, therefore,

$$1+\sum _{i=1}^{k} w^{i^2}=-\sum _{i=1}^{k} w^{-i^2}. $$

\noindent So,

$$\arrowvert 1+\sum _{i=1}^{k} w^{i^2} \arrowvert = \arrowvert \sum _{i=1}^{k} w^{-i^2} \arrowvert = \arrowvert \sum _{i=1}^{k} w^{i^2}  \arrowvert. $$\\

\noindent Also, since $n$ is prime, therefore, for each non-zero $j \in \mathbb{Z}_n$ 

$$1+\sum _{i=1}^{k} w^{ji^2}=-\sum _{i=1}^{k} w^{-ji^2}. $$

\noindent Therefore,

$$\arrowvert 1+\sum _{i=1}^{k} w^{ji^2} \arrowvert = \arrowvert \sum _{i=1}^{k} w^{-ji^2} \arrowvert = \arrowvert \sum _{i=1}^{k} w^{ji^2}  \arrowvert. $$

\noindent Since $\arrowvert 1+z \arrowvert = \arrowvert z \arrowvert $ if and only if $Re(z)=-\dfrac{1}{2} $, therefore, we get

$$Re(w^{1^{2}}+w^{2^{2}}+ \ldots +w^{k^{2}})=-\dfrac{1}{2}.$$

\noindent Then, $Re(1+w^{1^{2}}+w^{2^{2}}+ \ldots +w^{k^{2}})=\dfrac{1}{2} $, such that

$$\arrowvert 1+\sum _{i=1}^{k} w^{i^2} \arrowvert = \arrowvert \sum _{i=1}^{k} w^{i^2} \arrowvert. $$

\noindent As $\lbrace f_{0}, f_{1},\ldots,f_{n-1} \rbrace$ is an equiangular cyclic $(n,k)$-frame, therefore for each non-zero $j \in \mathbb{Z}_n$,

$$\arrowvert \sum _{i=1}^{k} w^{i^2} \arrowvert = \arrowvert \sum _{i=1}^{k} w^{ji^2} \arrowvert.$$

\noindent Thus for each non-zero $j \in \mathbb{Z}_n$,

$$\arrowvert 1+ \sum _{i=1}^{k} w^{i^2} \arrowvert = \arrowvert 1+ \sum _{i=1}^{k} w^{ji^2} \arrowvert.$$

\noindent Hence

$$g_{1}=\dfrac{1}{\sqrt{n}}\left[
\begin{array}{c}
1\\
w^{1^2}\\ 
w^{2^2}\\
\vdots \\
w^{k^2}\\
\end{array}
\right]
$$\\

\noindent generates an equiangular cyclic $(n,k+1)$-frame.\\
\end{proof}

\begin{cor}
Let $n$ be a prime integer such that $n=2k+1$ where $k$ is odd. Then the set of quadratic residues form a difference set. Moreover, the set of residues together with $\lbrace 0 \rbrace $ also form a difference set.
\end{cor}

\noindent Let $n$ be a prime integer such that $n=2k+1$ where $k$ is odd. We have now seen two different ways of generating equiangular cyclic $(n,k+1)$-frames. Let us now look at these two frames closely. \\

\noindent Let $G$ be the frame generated by 

$$g_{1}=\dfrac{1}{\sqrt{n}}\left[
\begin{array}{c}
1\\
w^{1^2}\\ 
w^{2^2}\\
\vdots \\
w^{k^2}\\
\end{array}
\right]
$$\\

\noindent And, let $H$ be the frame generated by 

$$h_{1}=\dfrac{1}{\sqrt{n}}\left[
\begin{array}{c}
1\\
w^{-1^2}\\ 
w^{-2^2}\\
\vdots \\
w^{-k^2}\\
\end{array}
\right]
$$\\

\noindent Note that 

$$ w^{(n-j){i^2}}= w^{-ji^2}= w^{j({n-i^2})}.$$

\noindent Thus $g_j = h_{n-j} $ for every $j$. Therefore, the two frames $G$ and $H$ are just the permutations of each other, and hence, are equivalent in the sense of \cite{2}. \\

\section{Gauss Sums and Equiangular cyclic frames}

\noindent We shall now look at various properties of equiangular cyclic frames that are derived from the concept of Gauss sums. We first start with the definition of a Gauss sum (Gauss period). This is also useful in gaining some specific information on the projection matrices of above constructed equingular cyclic frames.\\

\begin{defn}
A Gauss sum is a sum of roots of unity written as 
$$\varphi(a,n)=\sum _{r \in \mathbb{Z}_n} e^{\frac{-i \pi {r^2} a}{n}} $$
\noindent where $a$ and $n$ are relatively prime integers.\\
\end{defn}

\noindent We look at the case $ a=-2$ and odd prime $n \equiv 3 \pmod 4$. These are also called {\em quadratic Gauss sums}. Then, 

$$\varphi(-2,n)=\sum _{r \in \mathbb{Z}_n} w^{r^2}. $$\\

\noindent Define $R$, $T$ and $N$ as 

$$R=\sum _{r=1}^{k} w^{r^2},\ \ \ \ \ \ T= \sum _{r=1}^{n-1}(\dfrac{r}{n})_{L} w^{r} \ \ \ \ \ \ \  and \ \ \ \ \ \ N=\sum _{r=1}^{k} w^{-r^2}.$$\\

\noindent Clearly, $N=\overline{R}$ and $\varphi(-2,n)= 1+2R $. \\

\noindent Also, $T=R-N$ and $1+R+N=0$.\\ 

\noindent Therefore, $T=N+R=1+2R=\varphi(-2,n).$ \\

\noindent Then Gauss showed that $ \varphi(-2,n) = \sqrt{n}i$. Let $R=a+bi$. So, 
$$
\begin{array}{rcl}
\sqrt{n}i & = & \varphi(-2,n)\\
 & = & 1+2R\\
 & = & (1+2a)+2bi 
\end{array}
$$

\noindent Let us now compare real and imaginary parts. We get $a=-\dfrac{1}{2}$ as shown earlier and $b=\dfrac{\sqrt{n}}{2}$. Hence, 

$$ R=-\dfrac{1}{2}+\dfrac{\sqrt{n}}{2}i, \ \ \ \ \ \ \ \ N=\overline{R}=-\dfrac{1}{2}-\dfrac{\sqrt{n}}{2}i $$ 

\noindent and,

$$|R|=|N|=\dfrac{\sqrt{n+1}}{2}.$$\\

\noindent Using Gauss sums we were able to find out the value of $R$ which is used in the following results.\\

\begin{prop}
Let $n$ be a prime integer such that $n=2k+1$ where $k$ is odd. Consider the equiangular cyclic $(n,k)$-frame generated by 

$$f_{1}=\dfrac{1}{\sqrt{n}}\left[
\begin{array}{c}
w^{1^2}\\ 
w^{2^2}\\
\vdots \\
w^{k^2}\\
\end{array}
\right]
$$

\noindent Then, 

$$
\langle f_{j},f_{i} \rangle = \left\{
\begin{array}{rcl}
\dfrac{1}{n}R & & (j-i)\ is\ a\ quadratic\ residue \\
 & & \\
\dfrac{1}{n}\overline{R} & & (j-i)\ is\ a\ quadratic\ nonresidue \\
\end{array}\right.
$$\\

\end{prop}

\begin{proof}
\noindent Let us consider the inner-product $\langle f_j, f_i \rangle $ as follows
$$
\langle {\dfrac{1}{\sqrt{n}}} \left[
\begin{array}{c}
w^{j1^2}\\ 
w^{j2^2}\\
\vdots \\
w^{jk^2}\\
\end{array}
\right],
{\dfrac{1}{\sqrt{n}}} \left[
\begin{array}{c}
w^{i1^2}\\ 
w^{i2^2}\\
\vdots \\
w^{ik^2}\\
\end{array}
\right] \rangle=\frac{1}{n} \sum_{r=0}^{k} w^{(j-i)r^2}
$$

\noindent Hence, we get
$$
\langle f_j, f_i \rangle = \left\{
\begin{array}{rcl}
\frac{1}{n}\sum_{r=0}^{k} w^{r^2} & & (j-i)\ is\ a\ quadratic\ residue \\ \\
\frac{1}{n}\overline{\sum_{r=0}^{k} w^{r^2}} & & (j-i)\ is\ a\ quadratic\ nonresidue \\ \\
\end{array}\right.
$$

\noindent Therefore, 
$$
\langle f_{j},f_{i} \rangle = \left\{
\begin{array}{rcl}
\frac{1}{n}R & & (j-i)\ is\ a\ quadratic\ residue \\ \\
\frac{1}{n}\overline{R} & & (j-i)\ is\ a\ quadratic\ nonresidue \\
\end{array}\right.
$$
\end{proof}

\noindent As these frames are equiangular, therefore we know that $\arrowvert \langle f_{j},f_{i} \rangle \arrowvert = c$ for all $i \neq j$, where $c$ is some constant. By \cite{2}, we know that for equiangular cyclic $(n,k)$-frames, the value of this constant is given by 

$$ c=\sqrt{\frac{k(n-k)}{n^2(n-1)}}.$$\\

\noindent Let us first consider the case when $(j-i)$ is a quadratic residue and let $\langle f_{j},f_{i} \rangle = c \lambda_{ij}$ such that $\arrowvert \lambda_{ij}\arrowvert=1$. Therefore, by the above proposition, we get 

$$ \sqrt{\frac{k(n-k)}{n^2(n-1)}}\lambda_{ij}= \frac{R}{n}= -\frac{1}{2n}+\frac{\sqrt{n}}{2n}i.$$\\

\noindent Now by comparing real and imaginary parts, we can show that 

$$ \lambda_{ij}= \frac{-1}{\sqrt{n+1}}+\sqrt{\frac{n}{n+1}}i .$$ \\

\noindent Then we can see that 

$$
\lambda_{ij} = \left\{
\begin{array}{rcl}
\frac{-1}{\sqrt{n+1}}+\sqrt{\frac{n}{n+1}}i & & (j-i)\ is\ a\ quadratic\ residue \\ 
\frac{-1}{\sqrt{n+1}}-\sqrt{\frac{n}{n+1}}i & & (j-i)\ is\ a\ quadratic\ nonresidue \\ 
\end{array}\right.
$$\\

\noindent Now, we use Gauss sums for prime $ n \equiv 1 \pmod 4 $ such that $ n=4k^2+1 $, where k is an odd integer. We know that for such prime numbers $n$, the set of quartic residues forms a difference set and hence, generates an equiangular cyclic $(n,k^2)$-frame.\\

\noindent Define
$$ Q_{S}=\sum _{r \in S} w^{r},$$
\noindent for some set $S$.\\

\noindent As shown before, the four orbits of $H_4$ form a partition of $\mathbb{Z}_n$. Let us denote the orbits as $\lbrace H_4, \overline{H}_4, N_4, \overline{N}_4 \rbrace$, where $N_4$ is the orbit $aH_4$ for some quadratic nonresidue $a$. As $H_4$ generates an equiangular cyclic frame, therefore, 

$$ |Q_{H_4}|=|Q_{\overline{H}_4}|=|Q_{N_4}|=|Q_{\overline{N}_4}|.$$\\

\noindent Also we know that $Q_{H_4}=\overline{Q_{\overline{H}_4}}$ and $Q_{N_4}=\overline{Q_{\overline{N}_4}}$. \\

\noindent Let $Q_{H_4}=a+bi$ and $Q_{N_4}=c+di$. Then, Gauss showed $1+4Q_{H_4}=\sqrt{n}+\sqrt{2n+2\sqrt{n}}$. \\

\noindent Now comparing the real and imaginary parts we get, $a=\frac{\sqrt{n}-1}{4}$ and $b=\frac{\sqrt{2n+2\sqrt{n}}}{4}$. \\

\noindent As shown above, $T=Q_{H_2}-Q_{\overline{H}_2}$. Then, $T=\sqrt{n}$. Since $Q_{H_2}=2a$ and  $Q_{\overline{H}_2}=2c$, therefore, $c=\frac{-\sqrt{n}-1}{4}$ and $d=\frac{\sqrt{2n-2\sqrt{n}}}{4}$.\\

\noindent By giving similar arguments as above, we can prove the following proposition.

\begin{prop}
Let $n$ be a prime integer such that $n=4k^2+1$ for some odd integer $k$. Consider the equiangular cyclic $(n,k^2)$-frame generated by 

$$f_{1}=\dfrac{1}{\sqrt{n}}\left[
\begin{array}{c}
w^{n_1}\\ 
w^{n_2}\\
\vdots \\
w^{n_{k^2}}\\
\end{array}
\right]
$$\\

\noindent where $\lbrace n_1, n_2, \dots, n_{k^2} \rbrace$ is the set of all quartic residues of $n$. Then, 

$$
\langle f_{j},f_{i} \rangle = \left\{
\begin{array}{lcl}
\frac{1}{n}Q_{H_4} & & if \ (j-i) \in {H_4} \\ \\
\frac{1}{n}\overline{Q_{H_4}} & & if \ (j-i) \in \overline{H}_4\\\\
\frac{1}{n}Q_{N_4} & & if \ (j-i) \in {N_4} \\\\
\frac{1}{n}\overline{Q_{N_4}} & & if \ (j-i) \in \overline{N}_4\\
\end{array}\right.\\
$$

\end{prop}

\noindent As these frames are also equiangular, therefore we know that $\langle f_{j},f_{i} \rangle = c$ for all $i \neq j$, where $c$ is some constant. In this case, 

$$ c=\sqrt{\frac{k^2 (n-k^2)}{n^2(n-1)}}.$$

\noindent Let us first consider the case when $(j-i)$ is in $H_4$  and let $\langle f_{j},f_{i} \rangle = c \lambda_{ij}$ such that $\arrowvert \lambda_{ij}\arrowvert=1$. Therefore, by proposition above, we get 
$$ \sqrt{\frac{k^2 (n-k^2)}{n^2(n-1)}}\lambda_{ij}= \frac{R}{n}= \frac{\sqrt{n}-1}{4n}+\frac{\sqrt{2n+2\sqrt{n}}}{4n}i.$$

\noindent Now by comparing real and imaginary parts, we can show that 
$$ \lambda_{ij}= \frac{\sqrt{n}-1}{\sqrt{3k^2 +1}}+\sqrt{\frac{2n+2\sqrt{n}}{3k^2 +1}}i .$$ 

\noindent Similarly, it can be shown that 

$$
\lambda_{ij} = \left\{
\begin{array}{lcl}
\frac{\sqrt{n}-1}{\sqrt{3k^2 +1}}+\sqrt{\frac{2n+2\sqrt{n}}{3k^2 +1}}i & & if \  (j-i) \in {H_4} \\ \\
\frac{\sqrt{n}-1}{\sqrt{3k^2 +1}}-\sqrt{\frac{2n+2\sqrt{n}}{3k^2 +1}}i & &  if \  (j-i) \in {\overline{H}_4} \\ \\
\frac{-\sqrt{n}-1}{\sqrt{3k^2 +1}}+\sqrt{\frac{2n-2\sqrt{n}}{3k^2 +1}}i & & if \  (j-i) \in {N_4} \\ \\
\frac{-\sqrt{n}-1}{\sqrt{3k^2 +1}}-\sqrt{\frac{2n-2\sqrt{n}}{3k^2 +1}}i & &  if \  (j-i) \in {\overline{N}_4} \\
\end{array}\right.
$$\\

\noindent We shall use this information in the following section.\\

\section{Random and Burst Errors}

\noindent Recall from \cite{2} that the operator norm of the $m \times m$ correlation matrix $( \langle f_{i_k}, f_{i_l} \rangle)^m _{k,l=1} $ gives the error of $m$-erasures occuring in locations $ \lbrace i_1, \dots, i_m \rbrace $.\\

\begin{defn}
A set of $m$-erasures is called a burst error if $ \lbrace i_1, i_2, \dots, i_m \rbrace $ are consecutive integers. For any arbitrary collection $\lbrace i_1, \dots, i_m \rbrace$, the set of $m$-erasures is called a random error.
\end{defn}

\noindent We shall now look at the characteristic properties of the correlation matrix $( \langle f_j, f_i \rangle)_{ij} $. In case of $j$ erasures, we consider the corresponding adjacent submatrix of order $j$. Recall from section $3$ that the correlation matrix of an equiangular cyclic $(n,k)$-frame is {\em Toeplitz}. Hence all principle submatrices with consecutive rows and columns are same. Therefore we get the following result. \\

\begin{prop}
For every $m$, the norm of the burst error for $m$-erasures is constant.
\end{prop} 

\noindent Since burst errors are a particular type of random errors, hence it seems natural to assume that the minimum random error would be very small as compared to the burst error. However, it need not be true. As can be seen by numerical computation, in case of equiangular cyclic $(11,5)$-frame generated by quadratic residues, the minimum random error $0.7106$ is less than the burst error $0.7611$. However in case of equiangular cyclic $(7,3)$-frame generated by quadratic residues, the minimum random error $0.7517$ is same as the burst error.\\ 

\noindent By \cite{2}, we know that the equiangular cyclic frames are all optimal for $2$-erasures. Let us now consider the case of $3$-erasures. \\

\noindent First consider the case when $n$ is prime such that $n=2k+1$, where $k$ is odd. Let us consider the case when the $\lbrace i, j, l \rbrace $ coefficients are lost. Then the $3 \times 3$ correlation matrix is given by 

$$
C_{i,j,l} =\left[
\begin{array}{ccc}
\langle f_i, f_i \rangle & \langle f_j, f_i \rangle & \langle f_l, f_i \rangle \\
\langle f_i, f_j \rangle & \langle f_j, f_j \rangle &  \langle f_l, f_j \rangle \\
\langle f_i, f_l \rangle & \langle f_j, f_l \rangle & \langle f_l, f_l \rangle \\
\end{array}
\right].
$$
\\

\noindent We shall now try to get some information about the norm of this correlation submatrix. We know that $ \langle f_a, f_b \rangle = \sum_{t=1}^{k} w^{(a-b)t^2}$. Therefore

$$
C_{i,j,l} =\left[
\begin{array}{ccc}
\frac{k}{n} & \sum_{t=1}^{k} w^{(j-i)t^2} & \sum_{t=1}^{k} w^{(l-i)t^2} \\
\sum_{t=1}^{k} w^{(i-j)t^2} & \frac{k}{n} & \sum_{t=1}^{k} w^{(l-j)t^2} \\
\sum_{t=1}^{k} w^{(i-l)t^2} & \sum_{t=1}^{k} w^{(j-l)t^2} & \frac{k}{n} \\
\end{array}
\right]
$$

\noindent It is shown in the above theorem that there are only two possible off-diagonal entries in the error matrix. Hence,

$$
C_{i,j,l} =\left[
\begin{array}{ccc}
\frac{k}{n} & \lambda_{ij} c & \lambda_{il} c \\
\bar{\lambda}_{ij} c & \frac{k}{n} &  \lambda_{jl} c \\
\bar{\lambda}_{il} c & \bar{\lambda}_{jl} c & \frac{k}{n}   \\
\end{array}
\right]
$$\\

\noindent Then we can write $C_{i,j,l} = \frac{k}{n} I+cJ $, where $J$ is given by

$$
J =\left[
\begin{array}{ccc}
0 & \lambda_{ij} & \lambda_{il} \\
\bar{\lambda}_{ij} & 0 &  \lambda_{jl}  \\
\bar{\lambda}_{il}  & \bar{\lambda}_{jl} & 0   \\
\end{array}
\right]
$$\\

\noindent where the only possible values for $\lambda_{ij}, \lambda_{il},\lambda_{jl} $ are $\lambda$ or $ \bar{\lambda} $ depending on whether $ (i-j),(i-l),(j-l) $ are quadratic residues or nonresidues respectively. Note that the eigenvalue of $C_{i,j,l}$ is given by $\frac{k}{n}+c\alpha$ where $\alpha$ is an eigenvalue of $J$. Computing the characteristic polynomial of $J$ we get 

$$-x^3+3x+2 Re(\lambda_{ij} \bar{\lambda}_{il} \lambda_{jl})=0. $$\\

\noindent It is easy to check that 

$$
 Re(\lambda_{ij} \bar{\lambda}_{il} \lambda_{jl}) = \left\{
\begin{array}{lcl}
Re(\lambda^3) & & \lambda_{ij} = \lambda_{jl} \neq \lambda_{il} \\\\
Re\lambda & & otherwise
\end{array}\right.
$$\\

\noindent Therefore, the matrix $C_{i,j,l}$ can have at most two distinct possible norms in this case. Also we get that following are the only inequivalent possible forms of matrices for $3$-erasures. These are obtained $D C_{i,j,l} D^{-1}$, where $D$ is a diagonal matrix chosen to make the off-diagonal entries of first row as $1$.

$$
(i)\left[
\begin{array}{ccc}
0 & 1 & 1 \\
1 & 0 &  \lambda  \\
1  & \bar{\lambda} & 0   \\
\end{array}
\right]
(ii)\left[
\begin{array}{ccc}
0 & 1 & 1 \\
1 & 0 &  \lambda^3  \\
1  & \bar{\lambda}^3 & 0   \\
\end{array}
\right]
$$

$$
(iii)\left[
\begin{array}{ccc}
0 & 1 & 1 \\
1 & 0 &  \bar{\lambda}  \\
1  & \lambda & 0   \\
\end{array}
\right]
(iv)\left[
\begin{array}{ccc}
0 & 1 & 1 \\
1 & 0 &  \bar{\lambda}^3  \\
1  & \lambda^3 & 0   \\
\end{array}
\right]
$$

\noindent Note that to obtain two different sets of eigenvalues, we only need to consider forms $(i)$ and $(ii)$. \\

\noindent As mentioned above, the only possible values for $\lambda_{ij}, \lambda_{il},\lambda_{jl} $ are $\lambda$ or $ \bar{\lambda} $ depending on whether $ (i-j),(i-l),(j-l) $ are quadratic residues or nonresidues respectively. We also know that for a burst error, the corresponding adjacent submatrix is {\em Toeplitz}. Let $j=i+1$ and $l=i+2$. Then we know that $\lambda_{ij}= \lambda_{jl}$.  As shown above, we get the following equivalent classes of $J$ in case of burst errors.

$$
(a)\left[
\begin{array}{ccc}
0 & \lambda & \lambda \\
\bar{\lambda} & 0 &  \lambda  \\
\bar{\lambda}   & \bar{\lambda} & 0   \\
\end{array}
\right]
(b)\left[
\begin{array}{ccc}
0 & \lambda & \bar{\lambda} \\
\bar{\lambda} & 0 &  \lambda  \\
\lambda & \bar{\lambda} & 0   \\
\end{array}
\right]
$$\\

\noindent Observe that to obtain form (a) or (b), we must have 2 as a quadratic residue or nonresidue respectively. Now by combining Theorem $2.5$ with this observation, we get the following result about the matrix $J$ defined above.

\begin{prop}
The matrix $J$ is of the form (a) or (b) whenever $n \equiv \pm 1 \pmod{8}$ or $n \equiv \pm 3 \pmod{8}$ respectively.\\
\end{prop}

\noindent Similarly the inequivalent forms of $E_{i_1, \dots , i_m}$ can be obtained for $m> 3$ and similar observations can be made.\\

\noindent Also, we can follow the same procedure to find such forms for $n \equiv 1 \pmod 4$ such that $n=4k^2 +1$, $k$ is odd.\\

\section{Inequivalent frames}

\noindent From section $2$, we know that there exist multiple types of difference sets. Thus, it is very natural to find out whether different types of difference sets generate equivalent or inequivalent equiangular cyclic $(n,k)$-frame. \\

\noindent From \cite{2}, we know that the errors for any $m$-erasures are the same for equivalent frames. Therefore, one way to find out is to check for the error of $3$-erasures in case more than one difference set exists for $(n,k)$. \\

\noindent There exist two distinct types of difference sets for $n=31,k=15$ and $n=43,k=21$ as seen in \cite{10}. On numerically computing the maximum of norms of all $3 \times 3$ correlation submatrices for equiangular cyclic $(31,15)$-frames, we see that the frames generated by type $H_6$ gives $0.6663$ and the one generated by type $Q$ gives $0.6555$. In case of $n=43,k=21$, the maximum of norms of all $3 \times 3$ correlation submatrices generated by type $H_6$ is $0.6426$ whereas for the frame generated by type $Q$ is $0.6321$. \\

\noindent As can be seen, in these cases, the frames generated by the set of quadratic residues is optimal for 3 erasures. The computations also show that the two frames generated by distinct difference sets need not be equivalent. Following is a list of cases from \cite{10}, where multiple difference sets exist. 

$$
\begin{array}{|c|c|c|c|}
\hline
n & k & Type & Difference Set\\
\hline
 & & & \\
31 & 15 & H_6 & 1,2,3,4,6,8,12,15,16,17,23,24,27,29,30\\
31 & 15 & Q & 1,2,4,5,7,8,9,10,14,16,18,19,20,25,28\\
 & & & \\
43 & 21 & H_6 & 1,2,3,4,5,8,11,12,16,19,20,21,22,27,32,33,35\\
 & & & 37,39,41,42\\
43 & 21 & Q & 1,4,6,9,10,11,13,14,15,16,17,21,23,24,25,31,35\\
 & & & 36,38,40,41\\
 & & & \\ 
 \hline
\end{array}
$$  \\

\section{Spherical 1-design}

\noindent \cite{2} states the characterization of spherical $1$-design and $2$-designs for real uniform frames. In the following theorem, we show that the same characterization holds for spherical $1$-designs in case of complex uniform frames. Let us begin by the definition of a spherical $1$-design.

\begin{defn}
A set of vectors $\lbrace v_1, v_2, \dots, v_n \rbrace$ forms a spherical $1$-design if and only if  $\lbrace v_1, v_2, \dots, v_n \rbrace$ is in the sphere of $\mathbb{C}^k$ and 
$$ \int f dS=\frac{1}{n}\sum_{i=1}^{n} f(v_i)$$
\noindent for all polynomials $f$ of degree $1$.\\
\end{defn}

\noindent Following is the characterization the spherical $1$-design for complex frames.\\

\begin{thm}
Let $\lbrace v_1, v_2, \dots, v_n \rbrace$ be a set of vectors in  $\mathbb{C}^k$. Then the set of vectors $\lbrace v_1, v_2, \dots, v_n \rbrace$ form a $1$-design if and only if 
$$\sum_{i=1}^{n} v_i =0. $$
\end{thm}

\begin{proof}
Let 
$$f(z)=a_0 +a_1 z + \dots + a_k z_k+b_1 \bar{z}_1+b_2 \bar{z}_2 + \dots + b_k \bar{z}_k  $$

\noindent be a polynomial of degree $1$. Let $z_i = x_i + y_i$. \\

\noindent Then $\bar{z}_i = x_i - y_i$. So,
$$ f(z)=a_0 +(a_1+b_1)x_1+\dots+(a_k +b_k)x_k+(a_1 - b_1)y_1 i+\dots+(a_k - b_k)y_k i. $$

\noindent Therefore, 
$$ \int f dS=a_0.$$

\noindent Let $a=(a_1 +b_1,(a_1 - b_1)i,\dots,a_k +b_k,(a_k - b_k)i)$. Let 

$$
v_j=\left[
\begin{array}{c}
c_{j1}+d_{j1}i\\
c_{j2}+d_{j2}i\\
\vdots\\
c_{jk}+d_{jk}i\\
\end{array}
\right]
$$

\noindent Then consider $v_{j} \rightarrow w_{j}$ such that 

$$
w_j=\left[
\begin{array}{c}
c_{j1}\\
d_{j1}\\
c_{j2}\\
d_{j2}\\
\vdots\\
c_{jk}\\
d_{jk}\\
\end{array}
\right]
$$

\noindent Then $f(v_i)=a_0+a.w_i$. Therefore, 

$$
\begin{array}{lcl}
\frac{1}{n}\sum_{i=1}^{n} f(v_i) & = & \frac{1}{n}\sum_{i=1}^{n}(a_0 +a.w_i)\\
 & = & a_0 +\frac{1}{n} a. \sum_{i=1}^{n} w_i \\
\end{array}
$$

\noindent Hence,
$$ \int f dS=a_0=\frac{1}{n}\sum_{i=1}^{n} f(v_i) $$ 

\noindent if and only if 

$$\sum_{i=1}^{n} v_i =0. $$
\end{proof}

\noindent Note that for a prime $n$, equiangular cyclic $(n,k)$-frames form a spherical $1$-design since 
$$\sum_{i=0}^n f_i =0. $$
\\\\\\

\bibliographystyle{amsplain}

\end{document}